\documentclass[submission]{FPSAC2020}

\usepackage{cite}
\usepackage{lipsum}
\usepackage[utf8]{inputenc}
\usepackage[polish, english]{babel}
\usepackage{tikz}
\usepackage{pgfplots}
\usepackage{todonotes}
\usepackage{subfiles}
\usepackage[labelfont={normalsize}]{caption,subfig}
\usepackage{amsmath,amssymb, amsthm}
\usepackage{graphicx}

\newcommand{\verker}
{
	\begin{tikzpicture}
	\begin{axis}[
		axis x line*=bottom,
  		axis y line*=left,
  		enlargelimits=false,
  		xmin=0,
		xmax=2.4,
		ymin=0,
		ymax=2.2
	]

	\addplot [
		domain=-2:2, 
		samples=1000, 
		color=black,
		style={ultra thick}
	]
	({1/pi*(x*rad(asin(x/2))+sqrt(4-x*x))+x/2},{1/pi*(x*rad(asin(x/2))+sqrt(4-x*x))-x/2});
	\addplot[
	only marks,
    color=black,
    mark=o,
	style={thick}
    ]
    coordinates {
    (0,2)
    };
    \addplot[
    color=black,
    mark=o,
	style={thick}
    ]
    coordinates {
    (2,0)
    };
	\fill[thick,draw=black,fill=white] ({100*(1/pi*(1.5*rad(asin(0.5/2))+sqrt(4-1.5*1.5))+0.5/2)},{100*(1/pi*(0.5*rad(asin(1.5/2))+sqrt(4-0.5*0.5))-0.5/2)}) circle (100*0.03) node[anchor=south] {\tiny \large \hspace{45pt}$G(0.6)$};
	\fill[thick,draw=black,fill=white] ({100*(1/pi*(0.5*rad(asin(1.5/2))+sqrt(4-0.5*0.5))-0.5/2)},{100*(1/pi*(1.5*rad(asin(0.5/2))+sqrt(4-1.5*1.5))+0.5/2)}) circle (100*0.03) node[anchor=south west] {\tiny \large \hspace{-10pt} $G(0.4)$};
	\fill[thick,draw=black,fill=white] ({100*(1/pi*(1*rad(asin(1/2))+sqrt(4-1*1))+1/2)},{100*(1/pi*(1*rad(asin(1/2))+sqrt(4-1*1))-1/2)}) circle (100*0.03) node[anchor=south] {\tiny \large \hspace{50pt}$G(0.8)$};
	\fill[thick,draw=black,fill=white] ({100*(1/pi*(1*rad(asin(1/2))+sqrt(4-1*1))-1/2)},{100*(1/pi*(1*rad(asin(1/2))+sqrt(4-1*1))+1/2)}) circle (100*0.03) node[anchor=south west] {\tiny \large \hspace{-8pt} $G(0.2)$};
	\fill[thick,draw=black,fill=white] ({100*(1/pi*(-2*rad(asin(-2/2))+sqrt(4-2*2))-2/2)},{100*(1/pi*(-2*rad(asin(-2/2))+sqrt(4-2*2))+2/2)}) circle (100*0.03) node[anchor=west] {\tiny \large \hspace{4pt}$G(0.0)$};
	\fill[thick,draw=black,fill=white] ({100*(1/pi*(2*rad(asin(2/2))+sqrt(4-2*2))+2/2)},{100*(1/pi*(2*rad(asin(2/2))+sqrt(4-2*2))-2/2)}) circle (100*0.03) node[anchor=south] {\tiny \large $G(1.0)$};

	\end{axis}
	\end{tikzpicture}
}

\newcommand{\trzy}
{
	\colorlet{CiemnoZielony}{green!50!black}
\subfloat[]{
\begin{tikzpicture}[scale=0.08]
\begin{scope}
\clip (0,0) rectangle (42.426407,84.852814);
\draw[red,ultra thick] plot[smooth,scale=14.142136] file {bumpingcurve.txt};
\end{scope}
\begin{scope} 
\clip[](35,0) -- (35,1) -- (30,1) -- (30,3) -- (27,3) -- (27,4) -- (25,4) -- (25,5) -- (23,5) -- (23,6) -- (21,6) -- (21,7) -- (19,7) -- (19,8) -- (17,8) -- (17,9) -- (15,9) -- (15,11) -- (14,11) -- (14,12) -- (13,12) -- (13,13) -- (12,13) -- (12,14) -- (11,14) -- (11,16) -- (9,16) -- (9,18) -- (7,18) -- (7,19) -- (6,19) -- (6,21) -- (5,21) -- (5,25) -- (4,25) -- (4,26) -- (3,26) -- (3,31) -- (2,31) -- (2,34) -- (0,34) -- (0,0);
\draw[black!30] (0,0) grid (35,34);
\end{scope}
\draw[](35,0) -- (35,1) -- (30,1) -- (30,3) -- (27,3) -- (27,4) -- (25,4) -- (25,5) -- (23,5) -- (23,6) -- (21,6) -- (21,7) -- (19,7) -- (19,8) -- (17,8) -- (17,9) -- (15,9) -- (15,11) -- (14,11) -- (14,12) -- (13,12) -- (13,13) -- (12,13) -- (12,14) -- (11,14) -- (11,16) -- (9,16) -- (9,18) -- (7,18) -- (7,19) -- (6,19) -- (6,21) -- (5,21) -- (5,25) -- (4,25) -- (4,26) -- (3,26) -- (3,31) -- (2,31) -- (2,34) -- (0,34) -- (0,0) -- cycle ;
\traj{27}

\draw(20,17) node {$T=1$}; 
\end{tikzpicture}
\label{traj:A}
}
\hfill
\subfloat[]{
\begin{tikzpicture}[scale=0.08]
\begin{scope}
\clip (0,0) rectangle (42.426407,84.852814);
\draw[red,ultra thick] plot[smooth,scale=14.142136] file {bumpingcurve.txt};
\end{scope}
\begin{scope}
\clip[](53,0) -- (53,1) -- (51,1) -- (51,2) -- (44,2) -- (44,3) -- (41,3) -- (41,4) -- (36,4) -- (36,5) -- (35,5) -- (35,6) -- (33,6) -- (33,7) -- (31,7) -- (31,8) -- (30,8) -- (30,9) -- (28,9) -- (28,10) -- (26,10) -- (26,12) -- (24,12) -- (24,13) -- (22,13) -- (22,15) -- (19,15) -- (19,17) -- (18,17) -- (18,18) -- (17,18) -- (17,20) -- (15,20) -- (15,22) -- (14,22) -- (14,24) -- (13,24) -- (13,26) -- (11,26) -- (11,28) -- (10,28) -- (10,29) -- (9,29) -- (9,32) -- (7,32) -- (7,33) -- (6,33) -- (6,36) -- (5,36) -- (5,39) -- (4,39) -- (4,41) -- (2,41) -- (2,47) -- (1,47) -- (1,52) -- (0,52) -- (0,0);
\draw[black!30] (0,0) grid (53,52);
\end{scope}
\draw[](53,0) -- (53,1) -- (51,1) -- (51,2) -- (44,2) -- (44,3) -- (41,3) -- (41,4) -- (36,4) -- (36,5) -- (35,5) -- (35,6) -- (33,6) -- (33,7) -- (31,7) -- (31,8) -- (30,8) -- (30,9) -- (28,9) -- (28,10) -- (26,10) -- (26,12) -- (24,12) -- (24,13) -- (22,13) -- (22,15) -- (19,15) -- (19,17) -- (18,17) -- (18,18) -- (17,18) -- (17,20) -- (15,20) -- (15,22) -- (14,22) -- (14,24) -- (13,24) -- (13,26) -- (11,26) -- (11,28) -- (10,28) -- (10,29) -- (9,29) -- (9,32) -- (7,32) -- (7,33) -- (6,33) -- (6,36) -- (5,36) -- (5,39) -- (4,39) -- (4,41) -- (2,41) -- (2,47) -- (1,47) -- (1,52) -- (0,52) -- (0,0) -- cycle ;
\traj{27,25,21,21,16,16,15,15,13,13,12,12,11,11}

\draw (30,20) node {$T=2$};
\end{tikzpicture}
\label{traj:B}
}
\hfill
\subfloat[]{
\begin{tikzpicture}[scale=0.08]
\begin{scope}
\clip (0,0) rectangle (42.426407,84.852814);
\draw[red,ultra thick] plot[smooth,scale=14.142136] file {bumpingcurve.txt};
\end{scope}
\begin{scope}
\clip[](80,0) -- (80,1) -- (73,1) -- (73,2) -- (65,2) -- (65,3) -- (63,3) -- (63,4) -- (61,4) -- (61,5) -- (58,5) -- (58,6) -- (53,6) -- (53,7) -- (50,7) -- (50,8) -- (49,8) -- (49,9) -- (47,9) -- (47,10) -- (44,10) -- (44,12) -- (42,12) -- (42,13) -- (39,13) -- (39,15) -- (38,15) -- (38,17) -- (36,17) -- (36,18) -- (33,18) -- (33,20) -- (29,20) -- (29,23) -- (28,23) -- (28,24) -- (27,24) -- (27,25) -- (25,25) -- (25,27) -- (23,27) -- (23,28) -- (22,28) -- (22,29) -- (21,29) -- (21,30) -- (20,30) -- (20,33) -- (18,33) -- (18,34) -- (17,34) -- (17,36) -- (16,36) -- (16,37) -- (15,37) -- (15,39) -- (14,39) -- (14,40) -- (13,40) -- (13,42) -- (12,42) -- (12,44) -- (9,44) -- (9,47) -- (8,47) -- (8,49) -- (7,49) -- (7,52) -- (6,52) -- (6,54) -- (5,54) -- (5,56) -- (4,56) -- (4,62) -- (3,62) -- (3,67) -- (2,67) -- (2,70) -- (1,70) -- (1,74) -- (0,74) -- (0,0);
\draw[black!30] (0,0) grid (80,74);
\end{scope}
\draw[](80,0) -- (80,1) -- (73,1) -- (73,2) -- (65,2) -- (65,3) -- (63,3) -- (63,4) -- (61,4) -- (61,5) -- (58,5) -- (58,6) -- (53,6) -- (53,7) -- (50,7) -- (50,8) -- (49,8) -- (49,9) -- (47,9) -- (47,10) -- (44,10) -- (44,12) -- (42,12) -- (42,13) -- (39,13) -- (39,15) -- (38,15) -- (38,17) -- (36,17) -- (36,18) -- (33,18) -- (33,20) -- (29,20) -- (29,23) -- (28,23) -- (28,24) -- (27,24) -- (27,25) -- (25,25) -- (25,27) -- (23,27) -- (23,28) -- (22,28) -- (22,29) -- (21,29) -- (21,30) -- (20,30) -- (20,33) -- (18,33) -- (18,34) -- (17,34) -- (17,36) -- (16,36) -- (16,37) -- (15,37) -- (15,39) -- (14,39) -- (14,40) -- (13,40) -- (13,42) -- (12,42) -- (12,44) -- (9,44) -- (9,47) -- (8,47) -- (8,49) -- (7,49) -- (7,52) -- (6,52) -- (6,54) -- (5,54) -- (5,56) -- (4,56) -- (4,62) -- (3,62) -- (3,67) -- (2,67) -- (2,70) -- (1,70) -- (1,74) -- (0,74) -- (0,0) -- cycle ;
\traj{27,25,21,21,16,16,15,15,13,13,12,12,11,11,11,11,9,9,9,9,9,9,9,9,8,8,8,8,8}

\draw (45,20) node {$T=3$}; 
\end{tikzpicture}
\label{traj:C}
}

}

\newcommand{\traj}[1]
{
	\foreach \y [count=\j] in {#1} 
	{
		\draw[fill=CiemnoZielony!30] (\y,\j-1) rectangle +(1,1);
	}
}

\newtheorem{theorem}{Theorem}
\DeclareMathOperator{\Pos}{Pos}
\DeclareMathOperator{\boxi}{box}
\newcommand{\floor}[1]{\lfloor #1 \rfloor}
\newcommand{\strzalka}{\hspace{-3pt}\uparrow}

\newcommand{\rysuj}
{	\begin{center}
		\subfloat[]{
			\begin{tikzpicture}[scale=0.9]
				\diagram{4,3,2}
			\end{tikzpicture}
			\label{fig:rysa}
		}
		\hspace{30pt}
		\subfloat[]{
			\begin{tikzpicture}[scale=0.9]
				\tableau{{1,2,4,7},{3,6,9},{5,8}}
			\end{tikzpicture}
		 \label{fig:rysb}	
		 }
	\end{center}
}

\newcommand{\steprsk}
{
\centering
\hfill
\subfloat[]{
	\begin{tikzpicture}[scale=0.75]
	\clip (-1.5,-0.5) rectangle (4.5,4.5);

	\draw (0,0) rectangle +(1,1); 
	\node at (0.5,0.5) {16}; 
	\draw (1,0) rectangle +(1,1); 
	\node at (1.5,0.5) {37}; 
	\draw (2,0) rectangle +(1,1); 
	\node at (2.5,0.5) {41}; 
	\draw (3,0) rectangle +(1,1); 
	\node at (3.5,0.5) {82}; 
	\draw (0,1) rectangle +(1,1); 
	\node at (0.5,1.5) {23}; 
	\draw (1,1) rectangle +(1,1); 
	\node at (1.5,1.5) {53}; 
	\draw (2,1) rectangle +(1,1); 
	\node at (2.5,1.5) {70}; 
	\draw (0,2) rectangle +(1,1); 
	\node at (0.5,2.5) {74}; 
	\draw (1,2) rectangle +(1,1); 
	\node at (1.5,2.5) {99};

	\end{tikzpicture}
	\label{fig:RSKa}
}
\hfill
\subfloat[]{
\begin{tikzpicture}[scale=0.75]
\clip (-1.5,-0.5) rectangle (4.5,4.5);
\fill[blue!10] (1,0) rectangle +(1,1);
\fill[blue!10] (1,1) rectangle +(1,1);
\fill[blue!10] (0,2) rectangle +(1,1);

\draw (0,0) rectangle +(1,1); 
\node at (0.5,0.5) {16}; 
\draw (1,0) rectangle +(1,1); 
\node at (1.5,0.5) {37}; 
\draw (2,0) rectangle +(1,1); 
\node at (2.5,0.5) {41}; 
\draw (3,0) rectangle +(1,1); 
\node at (3.5,0.5) {82}; 
\draw (0,1) rectangle +(1,1); 
\node at (0.5,1.5) {23}; 
\draw (1,1) rectangle +(1,1); 
\node at (1.5,1.5) {53}; 
\draw (2,1) rectangle +(1,1); 
\node at (2.5,1.5) {70}; 
\draw (0,2) rectangle +(1,1); 
\node at (0.5,2.5) {74}; 
\draw (1,2) rectangle +(1,1); 
\node at (1.5,2.5) {99};

\draw[->] (-0.3,-0.45) to[bend left=60] (-0.3,0.45);
\draw[->] (0,1) +(-0.3,-0.45) to[bend left=60] +(-0.3,0.45);
\draw[->] (0,2) +(-0.3,-0.45) to[bend left=60] +(-0.3,0.45);
\draw[->] (0,3) +(-0.3,-0.45) to[bend left=60] +(-0.3,0.45);

\tiny
\node[] at (-1,0) {18};
\node[] at (-1,1) {37};
\node[] at (-1,2) {53};
\node[] at (-1,3) {74};

\end{tikzpicture}
\label{fig:RSKb}
}
\hfill
\subfloat[]{
	\begin{tikzpicture}[scale=0.75]
	\clip (-1.5,-0.5) rectangle (4.5,4.5);
	\fill[blue!10] (1,0) rectangle +(1,1);
	\fill[blue!10] (1,1) rectangle +(1,1);
	\fill[blue!10] (0,2) rectangle +(1,1);
	\fill[blue!10] (0,3) rectangle +(1,1);
	
	\draw (0,0) rectangle +(1,1); 
	\node at (0.5,0.5) {16}; 
	\draw (1,0) rectangle +(1,1); 
	\node at (1.5,0.5) {18}; 
	\draw (2,0) rectangle +(1,1); 
	\node at (2.5,0.5) {41}; 
	\draw (3,0) rectangle +(1,1); 
	\node at (3.5,0.5) {82}; 
	\draw (0,1) rectangle +(1,1); 
	\node at (0.5,1.5) {23}; 
	\draw (1,1) rectangle +(1,1); 
	\node at (1.5,1.5) {37}; 
	\draw (2,1) rectangle +(1,1); 
	\node at (2.5,1.5) {70}; 
	\draw (0,2) rectangle +(1,1); 
	\node at (0.5,2.5) {53}; 
	\draw (1,2) rectangle +(1,1); 
	\node at (1.5,2.5) {99}; 
	\draw (0,3) rectangle +(1,1); 
	\node at (0.5,3.5) {74};

	\draw[->] (-0.3,-0.45) to[bend left=60] (-0.3,0.45);
	\draw[->] (0,1) +(-0.3,-0.45) to[bend left=60] +(-0.3,0.45);
	\draw[->] (0,2) +(-0.3,-0.45) to[bend left=60] +(-0.3,0.45);
	\draw[->] (0,3) +(-0.3,-0.45) to[bend left=60] +(-0.3,0.45);
	
	\tiny
	\node[] at (-1,0) {18};
	\node[] at (-1,1) {37};
	\node[] at (-1,2) {53};
	\node[] at (-1,3) {74};
	
	\end{tikzpicture}
	\label{fig:RSKc}
}

}

\newcommand{\tableau}[1]
{
	\foreach \x [count=\i] in {#1} 
	{
		\foreach \y [count=\j] in \x 
		{
			\draw[ultra thick] (\j+1,0) -- (0,0) -- (0,\i+0.5);
			\draw (\j-1,\i-1) rectangle +(1,1); 
			\node[] at (\j-0.5,\i-0.5) {\y};
		}	
	}
}

\newcommand{\diagram}[1]
{
	\foreach \x [count=\i] in {#1} 
	{
		\draw[ultra thick] (\x+1,0) -- (0,0) -- (0,\i+0.5);
		\foreach \j in {1,...,\x}
		\draw (\j-1,\i-1) rectangle +(1,1); 
	}
}

\title[Hydrodynamic limit of RSK]{Hydrodynamic limit \\ of the Robinson--Schensted--Knuth algorithm}

\author[Mikołaj Marciniak]{Mikołaj Marciniak\thanks{\href{mailto:marciniak@mat.umk.pl}{marciniak@mat.umk.pl}. Mikołaj Marciniak was partially supported by Narodowe Centrum Nauki, grant number 2017/26/A/ST1/00189 and Narodowe Centrum Badań i Rozwoju, grant number POWR.03.05.00-00-Z302/17-00.}\addressmark{1}}
\address{\addressmark{1}Interdisciplinary Doctoral School “Academia Copernicana”, Faculty of Mathematics and Computer Science, Nicolaus Copernicus University in Toruń, ul.~Chopina 12/18, 87-100 Toruń, Poland}

\received{\today}

\abstract{We investigate the evolution in time of the position of a fixed number in the insertion tableau when the Robinson--Schensted--Knuth algorithm is applied to a sequence of random numbers. When the length of the sequence tends to infinity, a~typical trajectory after scaling converges uniformly in probability to some deterministic curve.}


\begin{document}


\maketitle

\section{Introduction}
The final published version of this manuscript can be found in \cite{Marciniak22}.

\subsection{Notations}
\begin{figure}[h]
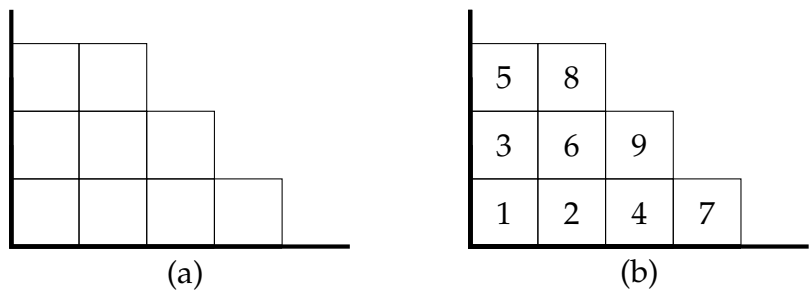

\rysuj
\caption{\protect\subref{fig:rysa} The Young diagram of shape $(4,3,2)\vdash 9$ and \protect\subref{fig:rysb} a standard Young tableau of shape $(4,3,2)\vdash 9$. }
\label{fig:diagtab}
\end{figure}

A \emph{partition} of a natural number $n$ is a break up of $n$ into a sum $n=\lambda_1+\lambda_2+\cdots+\lambda_k$ where $\lambda_1\geq \lambda_2\geq\cdots\geq\lambda_k>0$ are positive integer numbers. The vector \linebreak[4]{$\lambda=(\lambda_1,\lambda_2,\ldots,\lambda_k)$} is usually used to denote a partition. Let $\lambda\vdash n$ denote that $\lambda$ is a partition of a number $n$. A \emph{Young diagram} $\lambda=(\lambda_1,\lambda_2,\ldots,\lambda_k)$ is a finite collection of boxes arranged in left-justified rows with the row length $\lambda_j$ of the $j$-th row. Thus the Young diagram $\lambda$ is a graphical interpretation of the partition $\lambda$. A \emph{Young tableau} is a Young diagram filled with numbers. If the entries strictly decrease along each column from top to bottom and weakly increase along each row from left to right, a tableau is called \emph{semistandard}. A \emph{standard Young tableau} is a semistandard Young tableau with $n$ boxes which contains all numbers $1, 2, \ldots, n$. \cref{fig:diagtab} shows examples of a Young diagram and of a standard Young tableau.

The \emph{Robinson--Schensted--Knuth algorithm RSK} is a bijective algorithm which takes a finite sequence of numbers of length $n$ as the input and returns a pair of Young tableaux $(P,Q)$  with the same shape $\lambda\vdash n$. The semistandard tableau $P$ is called an \emph{insertion tableau}, and the standard tableau $Q$ is called a \emph{recording tableau}. In particular, the RSK algorithm assigns to any permutation $\sigma$ a pair $(P,Q)$ of standard Young tableaux. A detailed description of the RSK algorithm can be found in \cite[Chapter 1.6]{Romik2015}.

\begin{figure}[h]
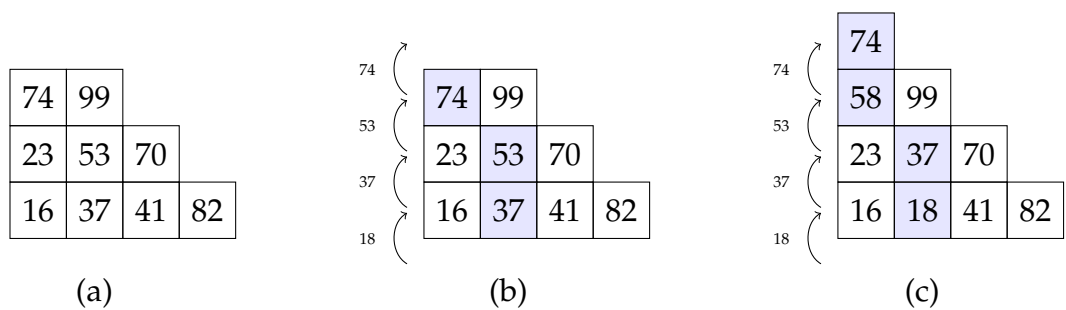

\steprsk
\caption{\protect\subref{fig:RSKa} The original tableau $P$.
	\protect\subref{fig:RSKb} The highlighted boxes form the bumping route which
	corresponds to an insertion of the number $18$. The numbers next to the arrows indicate the bumped numbers.\protect\subref{fig:RSKc} The
	output of the RSK insertion step. }
\label{fig:RSK}
\end{figure}

The RSK algorithm is based on applying the \emph{insertion step} to successive numbers from a given finite sequence $\{X_j\}_{j=1}^n$. The insertion step takes as input the previously obtained tableau $P(X_1, X_2, \ldots, X_{j-1})$ and the next number $X_j$ from the sequence. It produces as the output a new tableau $P(X_1, X_2, \ldots, X_j)$ with shape ``increased'' by one box; this tableau is obtained in the following way, see \cref{fig:RSK}. The RSK-insertion step starts in the first row with the number $x:=X_j$. The insertion step consists of inserting the number $x$ into the leftmost box in this row containing a number $y$ greater than $x$. We then move to the next
row with the number $y$ and repeat the action. At some row we are forced to insert the number at the end of the row, which will end the insertion step. The collection of rearranged boxes is called the \emph{bumping route}. 

We can say that the boxes with numbers are moved along the bumping route during an RSK insertion step. 

\subsection{Motivations}
The RSK algorithm is an important tool in algebraic combinatorics, especially in the context of Littlewood--Richardson coefficients and the plactic monoid\cite[Lecture 4]{Fulton1991}. 

For many years mathematicians have been studying the asymptotic behavior of the insertion tableau when we apply the Robinson--Schensted--Knuth algorithm to a random input. In the following paragraphs we will see several examples of such considerations. 

The Ulam--Hammersley problem \cite[Chapter 1.1]{Romik2015} concerns the typical length of the longest increasing subsequence in a random permutation. This problem corresponds to the problem of finding the typical length of the first row in the Young tableau obtained by the RSK algorithm from the sequence of independent random variables $\{X_j\}_{j=1}^n$ with the uniform distribution $U(0,1)$ on the unit interval $(0,1)$.

More generally, the RSK algorithm applied to the sequence of independent and identically distributed random variables with the uniform distribution $U(0,1)$ on the unit interval $(0,1)$ generates the \emph{Plancherel measure} on Young diagrams \cite[Chapter 1.8]{Romik2015}. The Plancherel measure is an important element of the representation theory of the symmetric groups because it describes how the left regular representation decomposes into irreducible components \cite[Chapter 3.3]{Fulton1991}.

Logan and Shepp \cite{LoganShepp1977} and Vershik and Kerov \cite{VershikKerov1986} described the limit shape of the insertion tableau $P(X_1, X_2, \ldots, X_n)$ obtained when we apply the RSK algorithm to a random finite sequence.

Romik and Śniady \cite{Romik2016} considered the limit shape of the bumping routes obtained when applying an RSK insertion step with a fixed number $w$ to an existing insertion tableau $P(X_1, X_2, \ldots, X_n)$ obtained from a random finite sequence. In \cite{Romik2015a} they also considered the limit shape of \emph{jeu de taquin} obtained from the recording tableau $Q(w, X_1, X_2, \ldots, X_n)$ made from a random finite sequence preceded by a fixed number $w$. 

\subsection{The main problem}
This paper concerns the asymptotic behavior of the insertion tableau when we apply the RSK algorithm to a random input. \emph{What can we say about the evolution over time of the insertion tableau from the viewpoint of box dynamics, when we apply the RSK algorithm to a sequence of independent random variables with the uniform distribution $U(0,1)$? How do the boxes move in the insertion tableau? If we investigate the scaled position of a box with a fixed number, will we get a deterministic limit, when the number of boxes tends to infinity?}

More specifically, we consider the insertion tableau $P(X_1, X_2, \ldots, X_n, w, X_{n+1}, \ldots, X_m)$ obtained by applying the
RSK algorithm to a random finite sequence containing a fixed number $w$ at some index. The box with this fixed number $w$ is being bumped by the RSK insertion step along the bumping routes. We will describe the scaled limit position of the box with the number $w$ depending on the ratio of the numbers $m$ and $n$. This problem was also stated by Duzhin \cite{Duzhin2019}.

Our main result for describing the scaled limit
position is \cref{thm:mainthm}. A proof of
\cref{thm:mainthm} is given in \cref{section:mainres}.

\section{Tools}
\subsection{The function $G(x)$}
\label{wzory}
Following Romik and Śniady in \cite[Chapter 5.1]{Romik2015a} we define the functions $F_{SC}, \Omega_{\star}, u, v$ needed to describe the macroscopic position of the new box in the RSK insertion step and additionally the function $G$
\begin{align*}
F_{SC}(y)&=\frac{1}{2}+\frac{1}{\pi}\bigg(\frac{y\sqrt{4-y^2}}{4}+\sin^{-1}\Big(\frac{y}{2}\Big)\bigg)\hspace{30pt}&\big(-2\leq y \leq 2\big),\\
\Omega_{\star}(y)&=\frac{2}{\pi}\bigg(\sqrt{4-y^2}+y\sin^{-1}\Big(\frac{y}{2}\Big)\bigg)&\big(-2\leq y \leq 2\big),\\
u(x)&=F_{SC}^{-1}(x) \hspace{30pt}&\big(0\leq x \leq 1\big),\\
v(x)&=\Omega_{\star}(u(x)) \hspace{30pt}&\big(0\leq x \leq 1\big),\\
G(x)&=\Big(\frac{v(x)+u(x)}{2},\frac{v(x)-u(x)}{2}\Big) \hspace{30pt}&\big(0\leq x \leq 1\big).
\end{align*}
The function $F_{SC}$ is the cumulative distribution function of \emph{Wigner's semicircle distribution}. The function $\Omega_{\star}$ is the limit shape of the Young tableau sampled from the \emph{Plancherel measure}. This curve is called the Logan--Shepp--Vershik--Kerov curve \cite{LoganShepp1977, VershikKerov1986}. The function $x \mapsto \big(u(x), v(x)\big)$ is a special parameterisation of the function $\Omega_{\star}$ and describes the limit position of the new box in the RSK insertion step applied with the number $x$ to the random Young tableau \cite[Chapter 5.1]{Romik2015a}. The continuous function $G\colon [0,1] \to \mathbb{R}_+^2$ is the function $\big(u(x), v(x)\big)$ rotated by 45 degrees (see \cite[Chapter 5.1]{Romik2015a}).
\begin{figure}
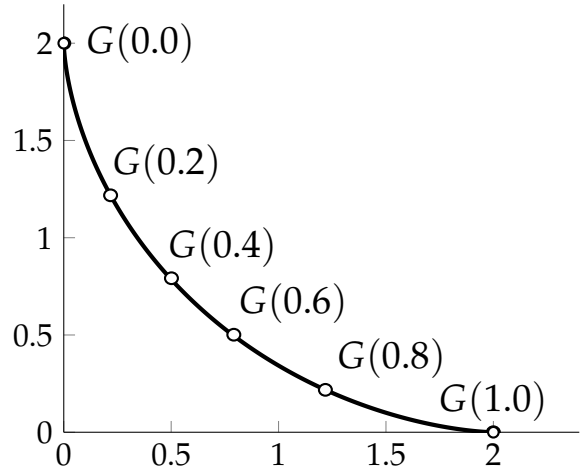

 \begin{center}
\verker
\end{center}
\caption{The graph of the function $G(x)$ with specified values for the argument\\
$x=0.0$, $x=0.2$, $x=0.4$, $x=0.6$ $x=0.8$ and $x=1.0$.}
\label{fig:gfunction}
\end{figure}
\cref{fig:gfunction} shows the graph of the function $G\colon [0,1] \to \mathbb{R}_+^2$. \\
\pagebreak[4] 
\subsection{The result of Romik and Śniady}
In the proof of \cref{thm:mainthm} we will need the following result of Romik and Śniady \cite[Theorem 5.1]{Romik2015a}. 

Let $\{X_j\}_{j=1}^\infty$ be a  sequence of independent random variables with the uniform distribution $U(0,1)$. Let $\square_n(x)\in \mathbb{N}^2$ denote the position of the box with the maximal number in the recording tableau $Q(X_1, \ldots, X_n, x)$ when we apply the RSK insertion step for the number $x\in[0,1]$ to the previously obtained tableau $P(X_1, \ldots, X_n)$. 
\begin{theorem}
\label{thm:romthm}
For each $x\in[0,1]$ the position $\square_n(x)$, after scaling by $\frac{1}{\sqrt{n}}$, converges in probability to a specific point $G(x)\in \mathbb{R}_+^2$, when $n$ tends to infinity: $$\frac{\square_n(x)}{\sqrt{n}}\overset{p}{\to}G(x).$$ 
\end{theorem}

\subsection{The partial order on the plane}
\label{order}
We define \emph{the partial order $\prec$ on the plane} as follows:
$(x_1,y_1)\prec(x_2,y_2)$ if and only if $x_1\leq x_2$ and $y_1\geq y_2$. 

For example the function $G$ is \emph{increasing with respect to the relation $\prec$}, i.e. if $x_1\leq x_2$ then $G(x_1)\prec G(x_2)$. Likewise, as
each row of the insertion tableau is non-decreasing, for each natural number $n\in \mathbb{N}$ the function $\square_n$ is increasing with respect to the relation $\prec$, i.e. if $x_1\leq x_2$ then $\square_n(x_1)\prec \square_n(x_2)$.
\subsection{The random increasing sequence}
\label{incseq}
Let $w\in(0,1]$ be a fixed number and let $\{X'_j\}_{j=1}^{m'}$ be a finite sequence of independent random variables with the uniform distribution $U(0,w)$. For $j\in\{1, 2, \ldots, m'\}$ we define $z(j)$ as the $j$-th order statistic, i.e. the $j$-th smallest number among $X'_1, X'_2, \ldots, X'_{m'}$. The sequence $\{z(j)\}_{j=1}^{m'}$ contains all elements of the sequence $\{X'_j\}_{j=1}^{m'}$ in the ascending order. The sequence $\{z(j)\}_{j=1}^{m'}$ will be called a \emph{random increasing sequence with the uniform distribution on the interval $[0,w]$}. 

In addition, the sequence $\{X'_j\}_{j=1}^{m'}$ is some permutation $\Pi=\left(\Pi_1, \Pi_2, \ldots, \Pi_{m'}\right)$ of the sequence $\{z(j)\}_{j=1}^{m'}$. Thus
$$\Big\{X'_j\Big\}_{j=1}^{m'}=\Big\{z(\Pi_j)\Big\}_{j=1}^{m'}=z\circ\Pi,$$
where $z$ is the function that acts pointwise on every element of the permutation $\Pi$. The permutation $\Pi$ is a random permutation with the uniform distribution. We see that a finite sequence of independent random variables with the uniform distribution $U(0,w)$ has the same distribution as a finite random increasing sequence with the uniform distribution on the interval $[0,w]$ permuted by a random permutation with the uniform distribution.

\section{The main result}
\label{section:mainres}
\subsection{Statement of Main Result}
\begin{figure}[h]
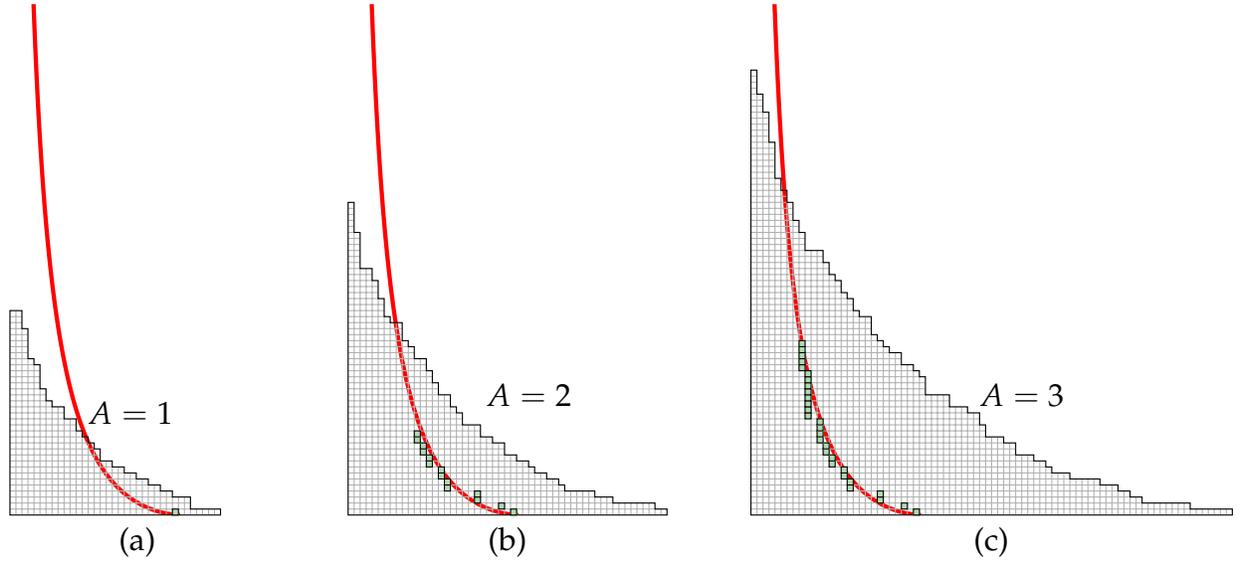

\trzy
\caption{ \protect\subref{traj:A}
The initial shape of the insertion tableau $P(X_1,\dots,X_n,w)$ immediately after the new box
with the number $w$ was added (the higlighted box in the bottom row) for $n=400$ and $w=0.5$.
\protect\subref{traj:B}
The shape of the insertion tableau 
$P(X_1,\dots,X_n,w,X_{n+1},\dots,X_{\lfloor T n \rfloor})$
at the time parameter $T=2$. The highlighted boxes indicate the trajectory of the box with the number $w$. The red smooth curve is the plot of $H$.
\protect\subref{traj:C} Analogous picture for $T=3$.
}
\label{fig:trzy}
\end{figure}

Our main result is \cref{thm:mainthm} describing the asymptotic behavior of the box with a fixed number. It states that when the number of boxes tends to infinity then the (scaled down) trajectory of the box converges in probability to the curve
$H\colon [1,\infty) \to \mathbb{R}_+^2$ given by
$$H(T):=\sqrt{T}\ G\left(\frac{1}{T}\right).$$ 
The same curve $H$ also happens to be the limit shape of the bumping routes \cite{Romik2016} in the RSK algorithm. \cref{fig:trzy} shows the graph of the curve $H$ and the experimentally determined trajectory of the box with the number $w=0.5$.

\pagebreak[3]

More specifically, let $w\in(0,1]$ be a fixed number. Let $\{X_j\}_{j=1}^\infty$ be a  sequence of independent random variables with the uniform distribution $U(0,1)$ on the unit interval $[0,1]$. For every $n\in\mathbb{N}$ we define the function $\Pos_n:\{n+1, n+2, \ldots\}\to\mathbb{N}^2$ by:
$$\Pos_n\left(j\right)=\boxi_w\biggl(P\left(X_1, \ldots, X_n, w, X_{n+1}, \ldots, X_j\right)\biggr)$$ for $j\in\{n+1, n+2, \ldots\}$, where for a tableau $P$ we denote by  $\boxi_w(P)\in\mathbb{N}^2$ the coordinates of the box with the number $w$.

\begin{theorem}
\label{thm:mainthm}
Let $R\in(1, \infty)$ be a real number. For each number $T\in[1,R]$ the random variable $\Pos_n(\floor{Tn})$, after scaling by $\frac{1}{\sqrt{wn}}$, converges in probability to $H(T)$, when $n$ tends to infinity. Moreover, the convergence is uniform, i.e.  for each $\epsilon>0$
$$	\lim_{n\to\infty}\mathbb{P}\Bigg(\sup_{T\in[1,R]} \left\|\frac{\Pos_n\left(\floor{Tn}\right)}{\sqrt{wn}}-H\left(T\right)\right\|>\epsilon\Bigg)=0.$$
\end{theorem}

\subsection{The strategy of the proof of \cref{thm:mainthm}}
The proof will be presented in consecutive subsections of this section. 
Here we present a sketch of the proof

First, we prove pointwise convergence in the following way.
\begin{itemize}
\item We reduce the case to numbers
less than $w$.
\item We restrict our attention only to the permutation generated by the numbers.  
\item We use the property of the RSK algorithm that the insertion tableau is equal to the recording tableau of the inverse permutation.
\item We return to the sequence of numbers from the interval $[0,1]$.
\item We apply the Romik-Sniady result about bumping routes to deduce the pointwise convergence.
 \end{itemize}
 
Second, using Dini’s second theorem we prove the uniform convergence.

\subsection{Pointwise convergence}

First, we prove only pointwise convergence, i.e. we prove that for each number $T\in[1,\infty)$ and for each $\epsilon>0$
$$	\lim_{n\to\infty}\mathbb{P}\Bigg( \left\|\frac{\Pos_n\left(\floor{Tn}\right)}{\sqrt{wn}}-H\left(T\right)\right\|>\epsilon\Bigg)=0.$$ 

\begin{proof}
We apply the RSK algorithm to a random sequence of real numbers containing the number $w$ and investigate the position of the box with the number $w$ in the insertion tableau. Any insertion step applied to a number greater that $w$ does not change the position of the number $w$ in the tableau, so it is enough to consider only the subsequence containing numbers not greater than $w$. 

Now we will use this observation in the proof. Let $m=\floor{Tn}$. The probability that the same number occurs twice in the sequence $(w, X_1, X_2, \ldots)$ is equal to 0, hence without losing generality we assume that the numbers $w, X_1, X_2, \ldots$ are all different. Let $\{X'_j\}_{j=1}^\infty$ be the subsequence of the sequence $\{X_j\}_{j=1}^\infty$ containing all elements of the sequence $\{X_j\}_{j=1}^\infty$ which are less than $w$. The sequence $\{X'_j\}_{j=1}^{\infty}$ is a sequence of independent random variables with the uniform distribution $U(0,w)$. 

Let $n'=n'(n)$ and $m'=m'(n)$ denote the number of elements, respectively, of the sequences $\{X_j\}_{j=1}^n$ and $\{X_j\}_{j=1}^m$ which are smaller than $w$. Then there is an equality:
\begin{align*}
\Pos_n\Big(\floor{Tn}\Big)&=\boxi_w\Big(P\left(X_1, \ldots, X_n, w, X_{n+1}, \ldots, X_{m}\right)\Big)\\&=\boxi_w\Big(P\left(X'_1, \ldots, X'_{n'}, w, X'_{n'+1}, \ldots, X'_{m'}\right)\Big).
\end{align*}
The random variable $n'=\sum_{j=1}^{n}[X_j<w]$ counts how many numbers from the sequence $\{X_j\}_{j=1}^{n}$ are smaller than $w$, so $n'$ is a random variable with the binomial distribution with parameters $n$ and $w$. We denote it $n'\sim B(n,w)$. Likewise, the random variable $m'-n'$ counts how many numbers from the sequence $\{X_j\}_{j=n+1}^{m}$ are smaller than $w$, so \linebreak[4]{$m'-n'\sim B(m-n,w)$.} Moreover, the random variables $n'$ and $m'-n'$ are independent, because the random variables $X_1, X_2, \ldots$ are independent. 

From the Strong Law of Large Numbers \cite[Theorem 2.4.1]{Durrett2019} we know that if $n$ tends to infinity then the following limits exist almost surely:
\begin{align*}
\lim_{n\to\infty}\frac{n'}{n}&=\mathbb{E}[X_j<w]=\mathbb{P}(X_j<w)= w, \\
\lim_{n\to\infty}\frac{m'-n'}{m-n}&=\mathbb{E}[X_j<w]=\mathbb{P}(X_j<w)=w. 
\end{align*}
Therefore, also the following limits exist almost surely
\begin{align*}
\lim_{n\to\infty}\frac{m'}{n'}
= 1+\lim_{n\to\infty}\frac{m'-n'}{n'}
= 1+\lim_{n\to\infty}\frac{m-n}{n}
=T
\end{align*}
and
\begin{align*}
\lim_{n\to\infty}\frac{m'}{n}=\lim_{n\to\infty}\frac{n'}{n}\frac{m'}{n'} =wT.
\end{align*}

We define the function $z:\{1, 2, \ldots m'\}\cup\{m'+1, n'+\frac{1}{2}\}\to [0,1]$ that assigns to a number $j\in\{1, 2, \ldots, m'\}$, the $j$-th smallest number among $X'_1, X'_2, \ldots, X'_{m'}$ and additionally $z\left(m'+1\right)=w$ and $z\left(n'+\frac{1}{2}\right)=\frac{z\left(n'\right)+z\left(n'+1\right)}{2}$. 

From the property of the random increasing sequence with the uniform distribution (\cref{incseq}) we have:
$$\Big\{X'_j\Big\}_{j=1}^{m'}=\Big\{z(\Pi_j)\Big\}_{j=1}^{m'}=z\circ\Pi,$$
where $\{z(j)\}_{j=1}^{m'}$ is a random increasing sequence with the uniform distribution on the interval $[0,w]$ and $\Pi$ is a random permutation of range $m'$ with the uniform distribution. The function $z$ acts pointwise on every element of the permutation $\Pi$. 
Similarly, the function $z$ acts on the insertion tableau by acting on each box individually. The function $z$ by acting on the input of the RSK algorithm changes pointwise the numbers in the insertion tableau, since the function $z$ is increasing. Then
\begin{align*}
\Pos_n\left(m\right)&=\boxi_w\Big(P\left(X'_1, \ldots, X'_{n'}, w, X'_{n'+1}, \ldots, X'_{m'}\right)\Big)
\\&= \boxi_w\bigg(P\Big(z(\Pi_1), \ldots, z\left(\Pi_{n'}\right), w, z\left(\Pi_{n'+1}\right), \ldots, z\left(\Pi_{m'}\right)\Big)\bigg)
\\&=\boxi_w\Big(z\circ P\left(\Pi_1, \ldots, \Pi_{n'}, m'+1, \Pi_{n'+1}, \ldots, \Pi_{m'}\right)\Big)
\\&=\boxi_{m'+1}\Big(P\left(\Pi_1, \ldots, \Pi_{n'}, m'+1, \Pi_{n'+1}, \ldots, \Pi_{m'}\right)\Big).
\end{align*}
We denote the inverse permutation to $\Pi$ by $\Pi^{-1}=\Big(\Pi^{-1}_1, \ldots, \Pi^{-1}_{m'}\Big)$ and we define the permutation $$\Pi\strzalka=\Big(\Pi_1, \ldots, \Pi_{n'}, m'+1, \Pi_{n'+1}, \ldots, \Pi_{m'}\Big)$$ as a natural extension of the permutation $\Pi$. Then $\Pi\strzalka^{-1}=\Big(\Pi\strzalka^{-1}_1, \ldots, \Pi\strzalka^{-1}_{m'+1}\Big)$ where
$$
\Pi\strzalka^{-1}_j=
\begin{cases}
\Pi_j^{-1}&\text{if}\hspace{10pt}\Pi_j^{-1}\leq n', \\
\Pi_j^{-1}+1& \text{if}\hspace{10pt}n'<\Pi_j^{-1}\leq m', \\
n'+1& \text{if}\hspace{10pt}j=m'+1. 
\end{cases}
$$ 
In addition, we will use the fact \cite{Schuetzenberger1963} that for any permutation $\Pi\strzalka$ the insertion tableau of $\Pi\strzalka$ is equal to the recording tableau of the inverse permutation $\Pi\strzalka^{-1}$:
$$P\left(\Pi\strzalka\right)=Q\left(\Pi\strzalka^{-1}\right).$$
Therefore
\begin{align*}
\Pos_n\left(m\right)&=\boxi_{m'+1}\bigg(P\left(\Pi_1, \ldots, \Pi_{n'}, m'+1, \Pi_{n'+1}, \ldots, \Pi_{m'}\right)\bigg)\\&=\boxi_{m'+1}\bigg(Q\left(\Pi\strzalka^{-1}\right)\bigg)\\&=\boxi_{m'+1}\bigg(Q\left(\Pi\strzalka^{-1}_1, \ldots, \Pi\strzalka^{-1}_{m'}, n'+1\right)\bigg).
\end{align*}
We want to rewrite this in terms of $\square_n(x)$ so that we can apply \cref{thm:romthm}. Now, by replacing the number $n'+1$ with the number $n'+\frac{1}{2}$ we get the inverse permutation $\Pi^{-1}$ 
\begin{align*}
\Pos_n\left(m\right)&=\boxi_{m'+1}\bigg(Q\Big(\Pi\strzalka^{-1}_1, \ldots, \Pi\strzalka^{-1}_{m'}, n'+1\Big)\bigg)\hspace{50pt}
\\&=\boxi_{m'+1}\bigg(Q\Big(\Pi\strzalka^{-1}_1, \ldots, \Pi\strzalka^{-1}_{m'}, n'+\frac{1}{2}\Big)\bigg)
\\&=\boxi_{m'}\bigg(Q\Big(\Pi^{-1}_1, \ldots, \Pi^{-1}_{m'}, n'+\frac{1}{2}\Big)\bigg).
\end{align*}

Since the function $z$ is increasing, when acting on the input
of the RSK algorithm it will not change the content of the recording tableau. Using the function $z$, we will get a sequence of independent random variables with the uniform distribution $U(0,1)$. Indeed
\pagebreak[2]
\begin{align*}
\Pos_n\left(m\right)&=\boxi_{m'}\bigg(Q\Big(\Pi^{-1}_1, \ldots, \Pi^{-1}_{m'}, n'+\frac{1}{2}\Big)\bigg)
\\&=\boxi_{m'}\bigg(z\circ Q\Big(\Pi^{-1}_1, \ldots, \Pi^{-1}_{m'}, n'+\frac{1}{2}\Big)\bigg)
\\&=\boxi_{m'}\Bigg(Q\bigg(z\Big(\Pi^{-1}_1\Big), \ldots, z\Big(\Pi^{-1}_{m'}\Big), z\Big(n'+\frac{1}{2}\Big)\bigg)\Bigg).
\end{align*}

The permutation $\Pi$ is a random permutation with the uniform distribution on $S_n$, so $\Pi^{-1}$ is also a random permutation with the uniform distribution on $S_n$. If we act with a random permutation on a random increasing sequence with the uniform distribution we will get a sequence of  independent random variables with the uniform distribution, thus the sequence $$z\circ\Pi^{-1}=\bigg(z\Big(\Pi^{-1}_1\Big), \ldots, z\Big(\Pi^{-1}_{m'}\Big)\bigg)$$ is a sequence of independent random variables with the uniform distribution $U(0,w)$. \\

We define the random variable $A_n$ and the sequence $\{Y_j\}_{j=1}^{m'}$: 
\begin{align*}
Y_j &=\frac{z\left(\Pi_j^{-1}\right)}{w} \qquad \text{for }j\in{1, 2, \ldots, m'},\\
A_n &=\frac{z\left(n'+\frac{1}{2}\right)}{w}=\frac{z(n')+z(n'+1)}{2w}.
\end{align*}
The sequence $\{Y_j\}_{j=1}^{m'}$ is a sequence of independent random variables with the uniform distribution $U(0,1)$, and the random variable $A_n$ converges almost surely to the limit $A=\frac{1}{T}$:
\begin{align*}
\lim_{n\to\infty}A_n&=\lim_{n\to\infty}\frac{z\left(n'\right)+z\left(n'+1\right)}{2w}\\&=\lim_{n\to\infty}\frac{z\left(n'\right)}{w}+\lim_{n\to\infty}\frac{z\left(n'+1\right)-z\left(n'\right)}{2w}\\&=\lim_{n\to\infty}\frac{n'}{m'}+0\\&=\frac{1}{T}.
\end{align*}
Therefore
\begin{align*}
\Pos_n\left(m\right)&=\boxi_{m'+1}\Bigg(Q\bigg(z\Big(\Pi^{-1}_1\Big), \ldots, z\Big(\Pi^{-1}_{m'}\Big), z\Big(n'+\frac{1}{2}\Big)\bigg)\Bigg)
\\&=\boxi_{m'+1}\Bigg(Q\bigg(\frac{z\left(\Pi^{-1}_1\right)}{w}, \ldots, \frac{z\left(\Pi^{-1}_{m'}\right)}{w}, \frac{z\left(n'+\frac{1}{2}\right)}{w}\bigg)\Bigg)
\\&=\boxi_{m'+1}\Big(Q\big(Y_1, \ldots, Y_{m'}, A_n\big)\Big)
\\&=\square_{m'}(A_n).
\end{align*}

Let $$G_{n}(x):=\frac{\square_{m'}(x)}{\sqrt{Tnw}}.$$ \\
The function $\square_{m'}$ is increasing with respect to the relation $\prec$ (of~\cref{order}).
Hence $G_{n}$ is also increasing with respect to the relation $\prec$, i.e. if $x_1\leq x_2$, then
\begin{align}
G_{n}(x_1)\prec G_{n}(x_2).
\label{ineq}
\end{align}
From \cref{thm:romthm} for each $x\in[0,1]$ the random variable $G_n(x)$ converges in probability to the limit $G(x)$, when $n$ tends to infinity. Indeed
\begin{align*}
G_{n}(x)=\sqrt{\frac{m'}{n}\frac{1}{Tw}} \hspace{5pt} \frac{\square_{m'}(x)}{\sqrt{m'}}\overset{p}{\to} \sqrt{\frac{Tw}{Tw}}\hspace{5pt}G(x)=G(x).
\end{align*}
\end{proof}

\subsection{The uniform convergence on the interval}

In this section we prove the uniform convergence on the interval $[1, R]$. 
\begin{proof}
For each $x$, the sequence $G_n(x)$ converges in probability to the limit $G(x)$, 
when $n$ tends to infinity. Both coordinates of each function $G_n(x)$ are monotonic 
and the limit function $G(x)$ is continuous. Then from Dini’s second theorem 
we get that for each $x$, the sequence $G_n(x)$ converges uniformly in probability 
to the limit $G(x)$, when $n$ tends to infinity.

Since $$\left\|G_n(A_n)-G(A)\right\| \leq \left\|G_n(A_n)-G(A_n)\right\| + \left\|G(A_n)-G(A)\right\|,$$ 
the sequence $G_n(A_n)$ converges uniformly in probability to the limit $G(A)$, when $n$ tends to infinity. 
Then for each $A\in[1,R]$ and each $\epsilon>0$ we have:
$$\lim_{n\to\infty}\mathbb{P}\bigg(\sup_{A\in[\frac{1}{R},1]}\left\|G_{n}(A_n)-G(A)\right\| 
> \epsilon\bigg)=0.$$
Using the inequality $T\leq R$ we obtain

$$\lim_{n\to\infty}\mathbb{P}\bigg(\sup_{A\in[\frac{1}{R},1]}\sqrt{T}\left\|G_{n}(A_n)-G(A)\right\| 
> \sqrt{R}\epsilon\bigg)=0,$$
$$	\lim_{n\to\infty}\mathbb{P}\Bigg(\sup_{T\in[1,R]} \left\|\frac{\Pos_n\left(\floor{Tn}\right)}
{\sqrt{wn}}-H\left(T\right)\right\|>\sqrt{R}\epsilon\Bigg)=0.$$
Thus for each $\epsilon>0$ we have
$$	\lim_{n\to\infty}\mathbb{P}\Bigg(\sup_{T\in[1,R]} \left\|\frac{\Pos_n\left(\floor{Tn}\right)}
{\sqrt{wn}}-H\left(T\right)\right\|>\epsilon\Bigg)=0.$$
\end{proof}

\subsection{Further questions}
Instead of one box with the number $w$ we can
simultaneously observe the position of $r$ boxes. 
If $r$ is finite, then \cref{thm:mainthm} describes
the hydrodynamic behavior of all $r$ boxes. Is this also true if $r$ goes to infinity? 

The function $H$ behaves asymptotically for $T\to\infty$ like the function $(\frac{1}{\sqrt{T}}, 2\sqrt{T})$. The function $H$ describes the asymptotic position of the box with a fixed number $w$. The $x$-coordinate of this position after scaling tends to zero. If we did not scale it, would it also be similar? Will the box with the number $w$ be moved to the first column? We expect that almost surely yes. How long will we have to wait for this? We expect that if the number $w$ appears in the RSK algorithm after $n$ steps, the waiting time for the number in the first column is $O(n^2)$. Is this really true? More formally, let $\{X_j\}_{j=1}^\infty$ be a sequence of independent random variables with the uniform distribution $U(0,1)$. What is the probability that the box with a fixed number $w$ will be in the first column in the insertion tableau $P(X_1, \ldots, X_n, w, X_{n+1}, \ldots, X_{\floor{Tn^2}})$ depending on parameter $T$?

\section*{Acknowledgments}
I would like to thank my supervisor, Professor Piotr Śniady, for his support, understanding, patience, positive outlook, guidance, valuable advice, and useful and constructive recommendations during writing this article. I would like to thank also to Dr Jacinta Torres and Dr Stephen Moore for the help with the text editing.

\bibliography{bib}{}
\bibliographystyle{amsalpha}

\end{document}